\documentclass[a4paper,10pt]{article}
\usepackage{stmaryrd}
\usepackage{amsfonts}
\usepackage{bbm}
\usepackage{amscd}
\usepackage{mathrsfs}
\usepackage{latexsym,amssymb,amsmath,amscd,amscd,amsthm,amsxtra,xypic}
\usepackage[dvips]{graphicx}
\usepackage[utf8]{inputenc}
\usepackage[T1]{fontenc}
\usepackage{lmodern}
\usepackage{amssymb}
\usepackage[all]{xy}
\usepackage{nicefrac,mathtools,enumitem}
\usepackage{microtype}
\usepackage{lmodern}

\newcommand {\emptycomment}[1]{}
\newtheorem{thm}{Theorem}[section]

\newtheorem{pro}[thm]{Proposition}
\newtheorem{ex}[thm]{Example}
\newtheorem{rmk}[thm]{Remark}
\newtheorem{defi}[thm]{Definition}

\newcommand{\lon }{\,\rightarrow\,}
\newcommand{\be }{\begin{equation}}
\newcommand{\ee }{\end{equation}}

\newcommand{\pf}{\noindent{\bf Proof.}\ }

\newcommand{\Real}{\mathbb R}

\newcommand{\CWM}{C^{\infty}(M)}

\newcommand{\frkg}{\mathfrak g}
\newcommand{\frkh}{\mathfrak h}

\newcommand{\frkl}{\mathfrak l}

\def\qed{\hfill ~\vrule height6pt width6pt depth0pt}

\newcommand{\br}[1]{   [ \cdot,    \cdot  ]   }

\newcommand{\idd}{\mathrm{id}}

\newcommand{\Ad}{\mathrm{Ad}}

\newcommand{\gl}{\mathfrak {gl}}

\begin{document}
\title{
{Hom-Lie groups of a class of Hom-Lie algebra
} }
\author{ Zhen Xiong \\
Department of Mathematics and Computer, Yichun University,\\
 Jiangxi, 336000, China
}
\date{}
\footnotetext{Supported by the NSF of China (No.11771382) and The Science and Technology Project(GJJ161029)of Department of Education, Jiangxi Province.}
\footnotetext{E-mail address:~205137@jxycu.edu.cn}

\maketitle
\textbf{Abstract}:  In this paper, the definition of Hom-Lie groups is given and one conntected component of Lie group $GL(V)$, which is not a subgroup of $GL(V)$, is a Hom-Lie group. More, we proved that there is a one-to-one relationship between Hom-Lie groups and Hom-Lie algebras $(\gl(V),[\cdot,\cdot]_\beta,\rm{Ad}_\beta)$.  Next, we also proved that if there is a Hom-Lie group homomorphism, then, there is a morphism between their Hom-Lie algebras. Last, as an application, we use these results on Toda lattice equation.
\textbf{MR(2010) Subject Classification}: 17B99, 22E99.


\section{Introduction}

The notion of Hom-Lie algebras was introduced by Hartwig, Larsson,
and Silvestrov in \cite{HLS} as part of a study of deformations of
the Witt and the Virasoro algebras.
Some $q$-deformations of the Witt and the Virasoro algebras have the
structure of a Hom-Lie algebra \cite{HLS,hu}. Because of close relation
to discrete and deformed vector fields and differential calculus
\cite{HLS,LD1,LD2}, more people pay special attention to this algebraic structure \cite{homlie1,MS2,yx,cy,cy2}.
Its geometric generalization is given in \cite{cls} and \cite{hom-Lie algebroids}.

In \cite{yx},the authors give a Hom-Lie algebra $(\gl(V),[\cdot,\cdot]_\beta,\rm{Ad}_\beta)$. This Hom-Lie algebra play an important role of studying structures of Hom-Lie algebras. Base on Hom-Lie algebra $(\gl(V),[\cdot,\cdot]_\beta,\rm{Ad}_\beta)$, there are many results are given: Omni-Hom-Lie algebra is given in \cite{yx}; Hom-big brackets are given in \cite{cy}; a Hom-Lie algebroid structure on $\varphi^!TM$ is given in \cite{cls}, and have the following results: there is a purely Hom-Poisson algebra structure on $\CWM$; a Hom-Lie algebra structure on the set of $(\sigma,\sigma)-$ derivations of an associative algebra is givn in \cite{cy2}, and so on.

As we know, from a Lie group, we get a Lie algebra; on the other hands, from a Lie algebra, we can have a Lie group. Hence, are there Hom-Lie groups, have similar results like Lie groups? In this paper, first, we study some properties of Hom-Lie algebra $(\gl(V),[\cdot,\cdot]_\beta,\rm{Ad}_\beta)$, then, give the definition of Hom-Lie groups, similar to Hom-Lie algebras, a Hom-Lie group twisted by a diffeomorphism is a Lie group. Next, on matrix space, we give some examples of Hom-Lie groups, and proved that one connected component of $GL(V)$, which is not a subgroup, is a Hom-Lie group. More, we proved that there is a one-to-one relationship a Hom-Lie group and a Hom-Lie algebra $(\gl(V),[\cdot,\cdot]_\beta,\rm{Ad}_\beta)$. Then, we give definitions of homomorphism on Hom-Lie groups, and proved that a homomorphism of Hom-Lie groups induce a morphism of Hom-Lie algebras. As an application, we study a deformation of Toda lattice equation.

The paper is organized as follow. In Section 2, we recall some
necessary background knowledge, including Hom-Lie algebras and morphism, Hom-Lie algebra $(\gl(V),[\cdot,\cdot]_\beta,\rm{Ad}_\beta)$. In Section 3, we study cohomology of Hom-Lie algebra $(\gl(V),[\cdot,\cdot]_\beta,\rm{Ad}_\beta)$, and have: cohomologies of Hom-Lie algebra $(\gl(V),[\cdot,\cdot]_\beta,\rm{Ad}_\beta)$ and Lie algebra $\gl(V)$ are isomorphic. In Section 4, we give definitions of Hom-Lie groups, homomorphism of Hom-Lie groups, and have mainly results: Theorem \ref{th1} and Theorem \ref{th2}. In  Section 5, for a class of deformations of Toda lattice equation, it is integrable.

\section{Preliminaries}

The notion of a Hom-Lie algebra was introduced in \cite{HLS}, see also \cite{MS2} for more information.
\begin{defi}
\begin{itemize}
\item[\rm(1.)]
  A Hom-Lie algebra is a triple $(\frkg,\br ,,\alpha)$ consisting of a
  vector space $\frkg$, a skewsymmetric bilinear map (bracket) $\br,:\wedge^2\frkg\longrightarrow
  \frkg$ and a linear transformation $\alpha:\frkg\lon\frkg$ satisfying $\alpha[x,y]=[\alpha(x),\alpha(y)]$, and the following hom-Jacobi
  identity:
  \begin{equation*}
   [\alpha(x),[y,z]]+[\alpha(y),[z,x]]+[\alpha(z),[x,y]]=0,\quad\forall
x,y,z\in\frkg.
  \end{equation*}

 A Hom-Lie algebra is called a regular Hom-Lie algebra if $\alpha$ is
a linear automorphism. When $\alpha=\idd$, Hom-Lie algebra $(\frkg,\br ,,\alpha)$ is just Lie algebra $(\frkg,\br ,)$.

 \item[\rm(2.)] A subspace $\frkh\subset\frkg$ is a Hom-Lie sub-algebra of $(\frkg,\br ,,\alpha)$ if
 $\alpha(\frkh)\subset\frkh$ and
  $\frkh$ is closed under the bracket operation $\br,$, i.e. for all $ x,y\in\frkh$,
  $[x,y] \in\frkh.  $
  \item[\rm(3.)] A morphism from the  Hom-Lie algebra
$(\frkg,[\cdot,\cdot]_{\frkg},\alpha)$ to the Hom-Lie algebra
$(\frkh,[\cdot,\cdot]_{\frkh},\delta)$ is a linear map
$\psi:\frkg\longrightarrow\frkh$ such that
$\psi([x,y]_{\frkg})=[\psi(x),\psi(y)]_{\frkh}$ and
$\psi\circ \alpha =\delta\circ \psi$.
  \end{itemize}
When $\psi$ is invertible, then $\psi$ is an isomorphism.
\end{defi}

The set of {\bf
$k$-cochains} on Hom-Lie algebra $(\frkg,\br ,,\alpha)$ with values in $V$, which we denote by
$C^k(\frkg;V)$, is the set of skewsymmetric $k$-linear maps from
$\frkg\times\cdots\times\frkg$($k$-times) to $V$:
$$C^k(\frkg;V):=\{\eta:\wedge^k\frkg\longrightarrow V ~\mbox{is a
linear map}\}.$$
Associated to the trivial representation, the set of $k$-cochains is $\wedge^k\frkg^*$. The corresponding coboundary operator $d:\wedge^k\frkg^*\longrightarrow\wedge^{k+1}\frkg^*$ is given by (see \cite{yx,xz})
$$d\xi(x_1,\cdots,x_{k+1})=\sum_{i<j}(-1)^{i+j}\xi([x_i,x_j],\alpha(x_1),\cdots,\hat{x_i},\cdots,\hat{x_j},\cdots,\alpha(x_{k+1})).$$

\begin{thm}\cite{yx}
Let $V$ be a vector space, and $\beta\in GL(V)$. Define a
skew-symmetric bilinear bracket operation
$[\cdot,\cdot]_\beta:\frkg\frkl(V)\times\frkg\frkl(V)\longrightarrow\frkg\frkl(V)$
by
$$[A,B]_\beta=\beta A\beta^{-1}B\beta^{-1}-\beta B\beta^{-1}A\beta^{-1},\quad \forall~ A,B\in\frkg\frkl(V),$$
where $\beta^{-1}$ is the inverse of $\beta$. Denote by
$\Ad_\beta:\frkg\frkl(V)\longrightarrow\frkg\frkl(V)$ the adjoint action on $\gl(V)$, i.e.
$\Ad_\beta(A)=\beta A\beta^{-1}.$
Then
$(\frkg\frkl(V),[\cdot,\cdot]_\beta,\Ad_\beta)$ is a regular Hom-Lie
algebra.
\end{thm}
Obviously, when $\beta=\idd$, $(\frkg\frkl(V),[\cdot,\cdot]_\beta,\Ad_\beta)$ is just Lie algebra $(\gl(V),[\cdot,\cdot])$.

\section{Properties of Hom-Lie algebra $(\frkg\frkl(V),[\cdot,\cdot]_\beta,\Ad_\beta)$ }
In this paper, we just study $\beta\circ\beta=\idd$, i.e. $\beta^{-1}=\beta$.
\begin{pro}
For any $C\in GL(V)$, let $\gamma=C\beta C^{-1}$, then there is an isomorphism from $(\frkg\frkl(V),[\cdot,\cdot]_\beta,\Ad_\beta)$ to $(\frkg\frkl(V),[\cdot,\cdot]_\gamma,\Ad_\gamma)$.
\end{pro}
\pf
Define map $F:(\frkg\frkl(V),[\cdot,\cdot]_\beta,\Ad_\beta)\longrightarrow (\frkg\frkl(V),[\cdot,\cdot]_\gamma,\Ad_\gamma)$ by $F(x)=CxC^{-1}$, then $F$ is an isomorphism.\qed
\begin{rmk}
When $\beta\neq\pm\idd$, there is not an isomorphism from $(\gl(V),[\cdot,\cdot])$ to $(\frkg\frkl(V),[\cdot,\cdot]_\beta,\Ad_\beta)$.
\end{rmk}
 For Hom-Lie algebra $(\gl(V),[\cdot,\cdot]_\beta,\Ad_\beta)$,
 coboundary operator $d:\wedge^k{\gl(V)}^*\longrightarrow\wedge^{k+1}{\gl(V)}^*$ is given by
$$
d\xi(x_1,\cdots,x_{k+1})=\sum_{i<j}(-1)^{i+j}\xi([x_i,x_j]_\beta,\Ad_\beta(x_1),\cdot,\widehat{x_{i,j}},\cdot,\Ad_\beta(x_{k+1})).
$$
For Lie algebra $(\gl(V),[\cdot,\cdot])$,  coboundary operator $\hat{d}:\wedge^k{\gl(V)}^*\longrightarrow\wedge^{k+1}{\gl(V)}^*$ is given by
$$\hat{d}\xi(x_1,\cdots,x_{k+1})=\sum_{i<j}(-1)^{i+j}\xi([x_i,x_j],x_1,\cdot,\widehat{x_{i,j}},\cdot,x_{k+1}).$$
Then, we have two chomology complexes: $(\oplus_k\wedge^k{\gl(V)}^*, d)$ of Hom-Lie algebra $(\gl(V),[\cdot,\cdot]_\beta,\Ad_\beta)$ and $(\oplus_k\wedge^k{\gl(V)}^*, \hat{d})$ of Lie algebra $(\gl(V),[\cdot,\cdot])$.
$\beta$  induce  map $\beta^r:\wedge^k{\gl(V)}^*\longrightarrow\wedge^{k}{\gl(V)}^*$ by
$$\beta^r(\xi)(x_1,\cdots,x_k)=\xi(x_1\beta,\cdots,x_k\beta).$$
And $\beta$ also induce map $\beta^l:\wedge^k{\gl(V)}^*\longrightarrow\wedge^{k}{\gl(V)}^*$ by
$$\beta^l(\xi)(x_1,\cdots,x_k)=\xi(\beta x_1,\cdots,\beta x_k).$$
\begin{pro}
For $\xi\in\wedge^k{\gl(V)}^*$, we have:
\begin{eqnarray}
\beta^r\circ d\xi&=&\hat{d}\beta^l(\xi);\label{3eq1}\\
d\beta^r(\xi)&=&\beta^l\circ\hat{d}\xi.\label{3eq2}
\end{eqnarray}
\end{pro}
\pf For $\xi\in\wedge^k{\gl(V)}^*$,
\begin{eqnarray*}
\beta^r\circ d\xi(x_1,\cdots,x_{k+1})&=&d\xi(x_1\beta,\cdots,x_{k+1}\beta)\\
&=&\sum_{i<j}(-1)^{i+j}\xi(\beta x_ix_j-\beta x_jx_i,\beta x_1,\cdot,\widehat{x_{i,j}},\cdot,\beta x_{k+1})\\
&=&\sum_{i<j}(-1)^{i+j}\beta^l(\xi)(x_ix_j-x_jx_i,x_1,\cdot,\widehat{x_{i,j}},\cdot,x_{k+1})\\
&=&\hat{d}\beta^l(\xi)(x_1,\cdots,x_{k+1})
\end{eqnarray*}
The proof of the rest of the conclusions is similar.\qed
\begin{rmk}In fact, we also have:
$\beta^l\circ \hat{d}\neq \hat{d}\circ \beta^l$, $\beta^r\circ d\neq d\circ \beta^r$.
\end{rmk}
Denote the set
of closed $k$-cochains of complex $(\oplus_k\wedge^k{\gl(V)}^*, d)$ by $Z^k(HL;\Real)$ and the set of exact $k$-cochains of complex $(\oplus_k\wedge^k{\gl(V)}^*, d)$ by $B^k(HL;\Real)$. Denote
the corresponding cohomology by
$$H^k(HL;\Real)=Z^k(HL;\Real)/B^k(HL;\Real).$$
Similarly, denote the set
of closed $k$-cochains of complex $(\oplus_k\wedge^k{\gl(V)}^*, \hat{d})$ by $Z^k(L;\Real)$ and the set of exact $k$-cochains of complex $(\oplus_k\wedge^k{\gl(V)}^*, \hat{d})$ by $B^k(L;\Real)$. Denote
the corresponding cohomology by
$$H^k(L;\Real)=Z^k(L;\Real)/B^k(L;\Real).$$

\begin{thm}
With the above notations, we have:
$$H^k(HL;\Real)=H^k(L;\Real).$$
\end{thm}
\pf For $\xi_1\in Z^k(HL;\Real)$, by $0=\beta^r\circ d\xi_1=\hat{d}\beta^l(\xi_1)$, we have: $\beta^l(\xi_1)\in Z^k(L;\Real)$. For $\xi_2\in Z^k(L;\Real)$,
by $\beta^r\circ d\beta^l(\xi_2)=\hat{d}\xi_2=0$, we have: $\beta^l(\xi_2)\in Z^k(HL;\Real)$. So, we have: $Z^k(HL;\Real)=Z^k(L;\Real).$

For $\eta_1\in B^k(HL;\Real)$, there is a $\eta_2\in \wedge^{k-1}{\gl(V)}^*$, and such that $d\eta_2=\eta_1$.
By $\hat{d}\beta^l(\eta_2)=\beta^r\circ d\eta_2=\beta^r(\eta_1)$, we have: $\beta^r(\eta_1)\in B^k(L;\Real)$.
On the other hand, for $\eta\in B^k(L;\Real)$, there is a $\eta_3\in \wedge^{k-1}{\gl(V)}^*$, and such that $\hat{d}\eta_3=\eta$.
By $d\beta^r(\eta_3)=\beta^l\hat{d}\eta_3=\beta^l(\eta)$, then, $\beta^l(\eta)\in B^k(HL;\Real)$. So, we have: $B^k(HL;\Real)=B^k(L;\Real)$.\qed

\section{Hom-Lie groups}
\begin{defi}
Let $G$ is a Lie group, $S$ is a submanifold of $G$. If there is a diffeomorphism $F:G\longrightarrow G$, such that $F(S)$ is a subgroup of $G$. Then $(S, F)$ is called a Hom-Lie group.  For $\forall x,y\in S$, if $F(y)$ is the inverse of $F(x)$ in Lie group $F(S)$, then we called $y$ is the Hom-inverse of $x$ in Hom-Lie group $(S, F)$. When $F=\idd$, Hom-Lie group $(S, F)$ is the Lie group $S$.
\end{defi}

\begin{ex}\label{ex1}
Matrix Lie group $GL(n;\mathbb{R})$, let $S_1=\{A\in GL(n;\mathbb{R})||A|<0\}$, $S_2=\{A\in GL(n;\mathbb{R})||A|>0\}$.  $S_1$ and $S_2$ are connected components of $GL(n;\mathbb{R})$, $S_2$ is a subgroup of $GL(n;\mathbb{R})$ and $S_1$ is not a subgroup of $GL(n;\mathbb{R})$.
We define the map $F:GL(n;\mathbb{R})\longrightarrow GL(n;\mathbb{R})$ by $F(A)=AP_{i,j}$, where $P_{i,j}$ is given by exchanging line $i$ and line $j$ of $\idd$, $P_{i,j}^2=\idd$. $F$ is a  diffeomorphism and $F^2=\idd$. So $(S_1, F)$ is a Hom-Lie group. $\forall A\in S_1$,  $P_{ij}A^{-1}P_{ij}$ is the Hom-inverse of $A$.
\end{ex}

\begin{ex}\label{ex2}
Matrix Lie group $O(n)$, let $S_3=\{A\in O(n)||A|=-1\}$, $S_4=\{A\in O(n)||A|=1\}$.  $S_3$ and $S_4$ are connected components of $O(n)$, $S_4$ is a subgroup of $O(n)$ and $S_3$ is not a subgroup of $O(n)$.
We define the map $F:O(n)\longrightarrow O(n)$ by $F(A)=AP_{i,j}$, where $P_{i,j}$ is given by exchanging line $i$ and line $j$ of $\idd$, $P_{i,j}^2=\idd$. $F$ is a diffeomorphism and $F^2=\idd$. So $(S_3, F)$ is a Hom-Lie group. $\forall A\in S_1$,  $P_{i,j}A^{-1}P_{i,j}$ is the Hom-inverse of $A$.
\end{ex}

\begin{defi}
A function $p:\mathbb{R}\longrightarrow (S,F)$ is called a one-parameter subgroup of Hom-Lie group $(S,F)$ if
\begin{itemize}
\item[1.]$p$ is continuous,
\item[2.]$F(p(0))=\idd$,
\item[3.]$F(p(t+s))=F(p(t))F(p(s))$.
\end{itemize}
\end{defi}

\begin{ex}\label{ex3}
For matrix Lie group $O(1;1)$, let
$$S_1=\{\left(
\begin{array}{clr}
cosh t&sinh t\\
sinh t&cosh t
\end{array}\right)|t\in\mathbb{R}\};
 S_2=\{\left(
\begin{array}{clr}
-cosh t&sinh t\\
sinh t&-cosh t
\end{array}\right)|t\in\mathbb{R}\};$$
$$S_3=\{\left(
\begin{array}{clr}
cosh t&-sinh t\\
sinh t&-cosh t
\end{array}\right)|t\in\mathbb{R}\};
S_4=\{\left(
\begin{array}{clr}
-cosh t&-sinh t\\
sinh t&cosh t
\end{array}\right)|t\in\mathbb{R}\}.$$
$S_1,S_2,S_3,S_4$ are connected components of $O(1;1)$.  $S_2,S_3,S_4$ are not subgroups of $O(1;1)$, $S_1$ is a subgroup of $O(1;1)$, hence $S_1$ is a Lie group.

We define the map $F_2:O(1;1)\longrightarrow O(1;1)$ by $F_2(A)=A(-\idd)$, then $(S_2,F_2)$ is a Hom-Lie group and $F_2^2=\idd$.

The map $F_3:O(1;1)\longrightarrow O(1;1)$ is given by
$F_3(A)=A\left(
\begin{array}{clr}
1&0\\
0&-1
\end{array}\right)$, then $(S_3,F_3)$ is a Hom-Lie group and $F_3^2=\idd$.

We define the map $F_4:O(1;1)\longrightarrow O(1;1)$ by
 $F_4(A)=A\left(
\begin{array}{clr}
-1&0\\
0&1
\end{array}\right)$, then $(S_4,F_4)$ is a Hom-Lie group and $F_4^2=\idd$.
Element of $S_i$ is a one-parameter subgroup of $(S_i,F_i), i=2,3,4$.
\end{ex}
$GL(V)$ is a matrix Lie group, $(\gl(V),[\cdot,\cdot])$ its Lie algebra, there is a map
$exp: \gl(V)\longrightarrow GL(V)$, for any $X\in \gl(V)$, $exp(X)=e^X$. More about the map $exp$, please see \cite{LG}.

For $\beta\in GL(V)$, let
$$M_{\beta}=\{e^{\beta X}\beta|X\in\gl(V)\},$$
then $M_{\beta}\subset GL(V)$, for $X\in\gl(V)$, let $p(t)=e^{t\beta X}\beta$, then $p(t)\subset M_{\beta}$.
\begin{pro}
 For $\beta\in GL(V),\beta^2=\idd$, let $M_{\beta}=\{e^{\beta X}\beta|X\in\gl(V)\}$, then $(M_{\beta},R_{\beta})$ is a Hom-Lie group.
\end{pro}
\pf When $\beta=\idd$, $M_{\beta}=GL(V)$ is a Lie group.

When $\beta\neq\idd$, if $X_0\in \gl(V)$ and such that $e^{\beta X_0}\beta=\idd$, we have $e^{\beta X_0}=\beta$, and $e^{2\beta X_0}=\beta^2=\idd$,
 then $X_0=0$, we have $\beta=\idd$, but $\beta\neq\idd$. So, $M_{\beta}$ is not a subgroup of $GL(V)$.

By $p(t)\subset M_{\beta}\subset\bigcup_{X\in\gl(V)} p(t)$, we have $M_{\beta}=\bigcup_{X\in\gl(V)} p(t)$. Because of $\beta\in\bigcap p(t)$, $M_{\beta}$ is a connected component of $GL(V)$. Then $M_{\beta}$ is a submanifold of $GL(V)$.

We define the map $R_{\beta}: GL(V)\longrightarrow GL(V)$ by $R_{\beta}(A)=A\beta$, $R_{\beta}$ is a diffeomorphism. Then $R_{\beta}\big(M_{\beta}\big)=\{e^{\beta X}|X\in\gl(V)\}$. So$(M_{\beta},R_{\beta})$ is a Hom-Lie group, and Hom-inverse of $e^{\beta X}\beta$ is $e^{-\beta X}\beta$.
\qed

For $X\in\gl(V)$, let $p(t)=e^{t\beta X}\beta$, then $p(t)$ is a one-parameter subgroup of Hom-Lie group $(M_{\beta},R_{\beta})$. For $Y\in\gl(V)$, we have
\begin{eqnarray}
\frac{d}{dt}(e^{t\beta X}\beta Ye^{-t\beta X}\beta)|_{t=0}&=&\beta X\beta Ye^0\beta-e^0\beta Y\beta X\beta \nonumber\\
&=&\beta X\beta Y\beta-\beta Y\beta XA\beta \nonumber\\
&=&[X,Y]_{\beta}.\label{eq1}
\end{eqnarray}
Actually, we proved the following results.
\begin{thm}\label{th1}
For $\beta\in GL(V), \beta^2=\idd$, there is a one-to-one relationship between  Hom-Lie group $(M_{\beta},R_{\beta})$ and regular Hom-Lie algebra $(\gl(V),[\cdot,\cdot]_{\beta},\rm{Ad}_{\beta})$. When $\beta=\idd$, Hom-Lie group $(M_{\beta},R_{\beta})$ is the Lie group $GL(V)$ and Hom-Lie algebra $(\gl(V),[\cdot,\cdot]_{\beta},\rm{Ad}_{\beta})$ is the Lie algebra $\gl(V)$.
\end{thm}

\begin{ex}\label{ext}
Let $M$ be a manifold and $\varphi:M\longrightarrow M$  a diffeomorphism. Define $\rm{Ad}_{\varphi^*}:\Gamma(\varphi^!TM)\longrightarrow\Gamma(\varphi^!TM)$ by
$$\rm{Ad}_{\varphi^*}(X)=\varphi^*\circ X\circ(\varphi^*)^{-1}, \forall X\in\Gamma(\varphi^!TM).$$
Define a skew-symmetric bilinear operation $[\cdot,\cdot]_{\varphi^*}:\wedge^2\Gamma(\varphi^!TM)\longrightarrow\Gamma(\varphi^!TM)
$ by
$$[X,Y]_{\varphi^*}=\varphi^*\circ X\circ(\varphi^*)^{-1}\circ Y\circ (\varphi^*)^{-1}-\varphi^*\circ Y\circ(\varphi^*)^{-1}\circ X\circ (\varphi^*)^{-1}.$$
Then, $(\varphi^!TM,[\cdot,\cdot]_{\varphi^*},\rm{Ad}_{\varphi^*},\idd)$ is a Hom-Lie algebroid and is called tangent Hom-Lie algebroid (\cite{cls}) .

In this example, now, we suppose $\varphi^2=\idd$ and $\varphi\neq\idd$, at a point $m\in M$, $x$ is a vector field on $M$,
if $\theta_t(m)$ is a one parameter Lie subgroup with respect to $x_m$, i.e. for any $f\in C^\infty(M)$, $\hat{d}f=x_mf=\lim_{t\rightarrow 0}\frac{1}{t}[f(\theta_t(m))-f(m)]$.

For tangent Hom-Lie algebroid, $X_m\in\Gamma(\varphi^!T_mM)$ and $X_m=x_{\varphi(m)}$, by representation of Hom-Lie algebroids on $C^\infty(M)$ (\cite{xz2}), for any $f\in C^\infty(M)$
\begin{eqnarray*}
d^1f(X_m)&=&(\varphi^*)^2\circ X_m\circ \varphi^*(f)\\
&=&X_m\circ(f\circ\varphi)\\
&=&x_{\varphi(m)}(f\circ\varphi)\\
&=&\lim_{t\rightarrow0}\frac{f\big(\varphi\circ\theta_t(\varphi(m))\big)-f(m)}{t}.
\end{eqnarray*}
Hence, $X_m$ with respect to one parameter Hom-Lie subgroup is $\varphi\circ\theta_t(\varphi(m))$.

On the other hands, $x_m\in T_mM$,
$$\frac{d}{dt}e^{tx_m}f|_{t=0}=x_mf.$$
$X_m\in\Gamma(\varphi^!T_mM)$, then $p(t)=e^{t\varphi^*\circ X_m}\varphi^*$ is one parameter Hom-Lie subgroup, and
\begin{eqnarray*}
\frac{d}{dt}e^{t\varphi^*\circ X_m}\varphi^*\circ f|_{t=0}&=&\varphi^*\circ X_m\circ\varphi^*\circ f\\
&=&x_{\varphi(m)}(f\circ\varphi).
\end{eqnarray*}

So, we got the same result through different paths.
\end{ex}
Actually, we can get the following results:
\begin{thm}
Let $M$ be a manifold, $\varphi:M\longrightarrow M$ is a diffeomorphism and $\varphi^2=\idd$, $\varphi\neq\idd$. If $x$ is a vector field on $M$, for any $m\in M$, then
$$D_m\varphi(x_m)\neq x_{\varphi(m)}.$$
\end{thm}
\pf
Let $\theta_t(m)$ is a one parameter Lie subgroup with respect to $x_m$. If $D_m\varphi(x_m)= x_{\varphi(m)}$, then we have
$$\varphi\circ\theta_t(m)=\theta_t(\varphi(m)),\quad i.e.,\quad \varphi\circ\theta_t(\varphi(m))=\theta_t(m).$$
But, by Example \ref{ext}, $\varphi\circ\theta_t(\varphi(m))$ is a one parameter Hom-Lie subgroup. Hence, $D_m\varphi(x_m)\neq x_{\varphi(m)}$.
\qed

\begin{defi}
Let $G_1$ and $G_2$ are Lie groups, $S_i$ is a submanifold of $G_i$, $i=1,2$. $(S_1,F_1)$ and $(S_2,F_2)$ are Hom-Lie groups. $\Phi:G_1\longrightarrow G_2$ is a Lie group homomorphism. We called $\Phi$ is a Hom-Lie group homomorphism from $(S_1,F_1)$ to $(S_2,F_2)$, if
\begin{itemize}
\item[1.)]$\Phi(S_1)\subset S_2$;
\item[2.)]$\forall x,y\in S_1$, $F_2\circ \Phi(xy)=\Phi\circ F_1(xy)$.
\end{itemize}
 \end{defi}
\begin{ex}
In Example \ref{ex1} and Example \ref{ex2}, we define map $\mathbb{I}:O(n)\longrightarrow GL(n;\mathbb{R})$ by $\mathbb{I}(A)=A$, i.e., $\mathbb{I}$ is an inclusion. $\mathbb{I}$ is a homomorphism from $(S_3,F)$ to $(S_1,F)$.
\end{ex}

\begin{ex}
The map $g:GL(n;\mathbb{R})\longrightarrow \mathbb{R}$ is given by $g(A)=|A|$, then $g$ is a Lie group homomorphism. In Example \ref{ex1}, $(S_1,F)$ is a Hom-Lie group, let $S_2=\{y\in \mathbb{R}|y<0\}$ and define the map $F_2:\mathbb{R}\longrightarrow\mathbb{R}$ by $F_2(y)=-y$, then $(S_2,F_2)$ is a Hom-Lie group. We have: $g(S_1)\subset S_2$ and $\forall A,B\in S_1, g\circ F(AB)=g(ABP_{i,j})=|A||B||P_{i,j}|=-|AB|=F_2\circ g(AB)$, so $g$ is a Hom-Lie group homomorphism.
\end{ex}

\begin{thm}\label{th2}
Let $\gamma\neq \beta, \gamma\neq \idd, \beta\neq\idd$, $\gamma^2=\idd=\beta^2$, $(M_{\gamma},R_{\gamma})$ and $(M_{\beta},R_{\beta})$ are Hom-Lie groups, with Hom-Lie algebras $(\gl(V),[\cdot,\cdot]_{\gamma},\rm{Ad}_{\gamma})$ and $(\gl(V),[\cdot,\cdot]_{\beta},\rm{Ad}_{\beta})$, respectively. Suppose that
$\Phi:(M_{\gamma},R_{\gamma})\longrightarrow(M_{\beta},R_{\beta})$ is a Hom-Lie group homomorphism. Then, there exists a real linear map $\phi:(\gl(V),[\cdot,\cdot]_{\gamma},\rm{Ad}_{\gamma})\longrightarrow(\gl(V),[\cdot,\cdot]_{\beta},\rm{Ad}_{\beta})$ such that
$$\Phi(e^{\gamma X})=e^{\beta\phi(X)}$$
for all $X\in\gl(V)$. The map $\phi$ has following additional properties:
\begin{itemize}
\item[(1.)]$\phi\circ\rm{Ad}_{\gamma}(X)=\rm{Ad}_{\beta}\circ\phi(X)$, for all $X\in\gl(V)$;
\item[(2.)]$\phi(e^{\gamma X}\gamma Ye^{-\gamma X}\gamma)=\Phi(e^{\gamma X}\gamma)\phi(Y)\Phi(e^{-\gamma X}\gamma)$, for all $X,Y\in\gl(V)$;
\item[(3.)]$\phi([X,Y]_{\gamma})=[\phi(X),\phi(Y)]_{\beta}$, for all $X,Y\in\gl(V)$;
\item[(4.)]$\beta\phi(X)=\frac{d}{dt}\Phi(e^{\gamma X})|_{t=0}$, for all $X\in\gl(V)$.
\end{itemize}
\end{thm}
\pf Since $\Phi$ is a Lie group homomorphism, then $R_{\beta}\circ\Phi(e^{\gamma X}\gamma)=\Phi\circ R_{\gamma}(e^{\gamma X}\gamma)$,
we have $\Phi(\gamma)=\beta$. And $\Phi$ is a continuous group homomorphism, for each $X\in\gl(V)$, so $\Phi(e^{t\gamma X})$ will be a one-parameter subgroup of Lie group $GL(V)$. Thus, there is a unique $Z\in\gl(V)$ such that
\begin{equation}\label{eq3}
\Phi(e^{t\gamma X})=e^{t\beta Z}.
\end{equation}
We define $\phi(X)=Z$ and check in several steps that $\phi$ has the required properties.

Step 1, $\Phi(e^{\gamma X})=e^{\beta\phi(X)}$.

This follows from (\ref{eq3}) and our definition of $\phi$, by putting $t=1$.

Step 2, $\phi(sX)=s\phi(X)$, for all $s\in\mathbb{R}$.

This is obviously.
Step 3, $\phi(X+Y)=\phi(X)+\phi(Y)$.

By Steps 1 and 2,
$$e^{t\beta\phi(X+Y)}=e^{\beta\phi(t(X+Y))}=\Phi(e^{t(\gamma X+\gamma Y)}).$$
By the Lie product formula and the fact that $\Phi$ is a continuous homomorphism, we have
\begin{eqnarray*}
e^t\beta\phi(X+Y)&=&\Phi\big(\lim_{m\rightarrow\infty}\big(e^t\gamma X/m e^t\gamma Y/m\big)^m\big)\\
&=&\lim_{m\rightarrow\infty}\big(\Phi(e^{t\gamma X/m})\Phi(e^{t\gamma Y/m})\big)^m.
\end{eqnarray*}
Then, we have
$$e^t\beta\phi(X+Y)=\lim_{m\rightarrow\infty}\big(e^{t\beta \phi(X)/m}e^{t\beta \phi(Y)/m}\big)^m=e^{t\beta(\phi(X)+\phi(Y))}.$$
Differentiating this result at $t=0$ gives the desired result.

Step 4, $\phi\circ\rm{Ad}_{\gamma}(X)=\rm{Ad}_{\beta}\circ\phi(X)$.

By Step 1,
\begin{eqnarray}
\Phi(e^{\gamma\rm{Ad}_{\gamma}(X)}\gamma)&=&\Phi(e^{\gamma\rm{Ad}_{\gamma}(X)})\Phi(\gamma)\nonumber\\
&=&e^{\beta\phi(\rm{Ad}_{\gamma}(X))}\beta.\label{eq4}
\end{eqnarray}
And $e^{\gamma\rm{Ad}_{\gamma}(X)}\gamma=e^{X\gamma}\gamma=\gamma e^{\gamma X}$, then
\begin{eqnarray}
\Phi(e^{\gamma\rm{Ad}_{\gamma}(X)}\gamma)&=&\Phi(\gamma e^{\gamma X})\nonumber\\
&=&\beta\Phi(e^{\gamma X})=\beta e^{\beta\phi(X)}.\label{eq5}
\end{eqnarray}
By (\ref{eq4}) and (\ref{eq5}), we have:
$$e^{\beta\phi(\rm{Ad}_{\gamma}(X))}=\beta e^{\beta\phi(X)}\beta=e^{\phi(X)\beta}.$$
Hence, $\phi(\rm{Ad}_{\gamma}(X))=\beta\phi(X)\beta=\rm{Ad}_{\beta}\circ\phi(X)$.

Step 5, $\phi(e^{\gamma X}\gamma Ye^{-\gamma X}\gamma)=\Phi(e^{\gamma X}\gamma)\phi(Y)\Phi(e^{-\gamma X}\gamma)$.
\begin{eqnarray*}
e^{\beta\phi\big(e^{t\gamma X}\gamma Ye^{-t\gamma X}\gamma \big)}&=&\Phi\big(e^{\gamma e^{t\gamma X}\gamma Ye^{-t\gamma X}\gamma}\big)\\
&=&\Phi\big(e^{ e^{tX\gamma} Y\gamma e^{-t X\gamma}}\big)\\
&=&\Phi\big(  e^{tX\gamma}\gamma\gamma e^{Y\gamma} e^{-tX\gamma}\big)\\
&=&\Phi\big( \gamma e^{t\gamma X}e^{\gamma Y}e^{-t\gamma X}\gamma\big)\\
&=&\beta\Phi(e^{t\gamma X})\Phi(e^{\gamma Y})\Phi(e^{-t\gamma X})\beta\\
&=&\beta\Phi(e^{t\gamma X})e^{\beta\phi(Y)}\Phi(e^{-t\gamma X})\beta\\
&=&e^{\beta\Phi(e^{t\gamma X})\beta\phi(Y)\Phi(e^{-t\gamma X})\beta}.
\end{eqnarray*}
Hence, we have
$$\phi\Big(e^{t\gamma X}\gamma Ye^{-t\gamma X}\gamma\Big)=\Phi(e^{t\gamma X}\gamma)\phi(Y)\Phi(e^{-t\gamma X}\gamma).$$
By putting $t=1$, we have the desired result.

Step 6, $\phi([X,Y]_{\gamma})=[\phi(X),\phi(Y)]_{\beta}$.

By Step 5, we have
\begin{eqnarray*}
\beta\phi\Big(e^{t\gamma X}\gamma Ye^{-t\gamma X}\gamma\Big)&=&\beta\Phi(e^{t\gamma X})\beta\phi(Y)\Phi(e^{-t\gamma X})\beta\\
&=&\beta e^{t\beta\phi(X)}\beta\phi(Y)e^{-t\beta\phi(X)}\beta.
\end{eqnarray*}
On the other hands,
\begin{eqnarray*}
\phi([X,Y]_\gamma)&=&\phi\Big(\frac{d}{dt}e^{t\gamma X}\gamma Ye^{-t\gamma X}\gamma|_{t=0}\Big)\\
&=&\frac{d}{dt}\phi(e^{t\gamma X}\gamma Ye^{-t\gamma X}\gamma)|_{t=0}\\
&=&\frac{d}{dt}e^{t\beta\phi(X)}\beta\phi(Y)e^{-t\beta\phi(X)}\beta|_{t=0}\\
&=&[\phi(X),\phi(Y)]_\beta.
\end{eqnarray*}

Step 6, $\beta\phi(X)=\frac{d}{dt}\Phi(e^{\gamma X})|_{t=0}$.

This follows our definition of $\phi$.\qed

\section{Applications}

 Consider the following equation
  \begin{equation}\label{eqe1}
  \frac{d}{dt}L=BL-LB=[B,L],
  \end{equation}
where $L$ is an $n\times n$ symmetric real tridiagonal matrix, and $B$ is
the skew symmetric matrix obtained from $L$ by
$$B=L_{>0}-L_{<0},$$
where $L_{>0(<0)}$ denotes the strictly upper (lower) triangular
part of $L$.  Based on Lie algebras,  Kodama,Y. and Ye, J.(\cite{KY,KY1}) study equation (\ref{eqe1}), and give an explicit formula for the solution to
the initial value problem.  Now we consider the following system:
\begin{eqnarray}
\frac{d}{dt}L &=& \beta B\beta L\beta-\beta L\beta B\beta \nonumber\\
 &=&[B,L]_\beta,\label{eqe2}
\end{eqnarray}
where $\beta^2=\idd$, then, from results of Section 4, we know that equation (\ref{eqe2}) is also integrable.


\end{document}